\documentclass[12pt]{article}
\usepackage{amssymb, amsmath, url, graphicx}

\def\3{\subset }
\def\4{\subseteq }
\def\<{\left<}
\def\>{\right>}

\def\n{\noindent }

\def\bit{\begin{itemize}}
\def\eit{\end{itemize}}
\def\3{\subset }
\def\4{\subseteq }

\def\calc{{\cal C}}
\def\cald{{\cal D}}
\def\calf{{\cal F}}
\def\0{\leqno}

\def\barr{\begin{array}}
\def\earr{\end{array}}
\def\dd{\displaystyle}

\def\Z{{\rlap{$\kern2pt{\rm Z}$}{\rm Z}\,}}

\title{\bf The number of chains of subgroups\\ of a finite elementary abelian $p$-group}
\author{Marius T\u arn\u auceanu}
\date{June 27, 2015}

\begin{document}

\maketitle

\begin{abstract}
In this short note we give a formula for the number of chains of
subgroups of a finite elementary abelian $p$-group. This completes
our previous work \cite{5}.
\end{abstract}



\section{Introduction}

Let $G$ be a group. A \textit{chain of subgroups} of $G$ is a set
of subgroups of $G$ totally ordered by set inclusion. A chain of
subgroups of $G$ is called \textit{rooted} (more exactly
$G$-\textit{rooted}) if it contains $G$. Otherwise, it is called
\textit{unrooted}. Notice that there is a bijection between the
set of $G$-\textit{rooted} chains of subgroups of $G$ and the set
of distinct fuzzy subgroups of $G$ (see e.g. \cite{5}), which is
used to solve many computational problems in fuzzy group theory.
\bigskip

The starting point for our discussion is given by the paper
\cite{5}, where a formula for the number of rooted chains of
subgroups of a finite cyclic group is obtained. This leads in
\cite{3} to precise expression of the well-known central Delannoy
numbers in an arbitrary dimension and has been simplified in
\cite{2}. Some steps in order to determine the number of rooted
chains of subgroups of a finite elementary abelian $p$-group are
also made in \cite{5}. Moreover, this counting problem has been
naturally extended to non-abelian groups in other works, such as
\cite{1,4}. The purpose of the current note is to improve the
results of \cite{5}, by indicating an explicit formula for the
number of rooted chains of subgroups of a finite elementary
abelian $p$-group.
\bigskip

Given a finite group $G$, we will denote by $\calc(G)$, $\cald(G)$
and $\calf(G)$ the collection of all chains of subgroups of $G$,
of unrooted chains of subgroups of $G$ and of $G$-rooted chains of
subgroups of $G$, respectively. Put $C(G)=|\calc(G)|$,
$D(G)=|\cald(G)|$ and $F(G)=|\calf(G)|$. The connections between
these numbers have been established in \cite{2}, namely:

\bigskip\n{\bf Theorem 1.} {\it Let $G$ be a finite group. Then
$$F(G)=D(G)+1 \mbox{ and }\, C(G)=F(G)+D(G)=2F(G)-1\,.$$}

In the following let $p$ be a prime, $n$ be a positive integer and
$\mathbb{Z}_p^n$ be an elementary abelian $p$-group of rank $n$
(that is, a direct product of $n$ copies of $\mathbb{Z}_p$). First
of all, we recall a well-known group theoretical result that gives
the number $a_{n,p}(k)$ of subgroups of order $p^k$ in
$\mathbb{Z}_p^n$, $k=0,1,...,n$.

\bigskip\n{\bf Theorem 2.} {\it For every $k=0,1,...,n$, we have
$$a_{n,p}(k)=\frac{(p^n-1)\cdots (p-1)}{(p^k-1)\cdots
(p-1)(p^{n-k}-1)\cdots (p-1)}\,.$$}Our main result is the
following.

\bigskip\n{\bf Theorem 3.} {\it The number of rooted chains of subgroups
of the elementary abelian $p$-group $\mathbb{Z}_p^n$ is
$$F(\mathbb{Z}_p^n)=2{+}2f(n)\sum_{k=1}^{n-1}\sum_{1\leq i_1<i_2<...<i_k\leq n-1}\frac{1}{f(n{-}i_k)f(i_k{-}i_{k-1})\cdots
f(i_2{-}i_1)f(i_1)}\,,$$where
$f:\mathbb{N}\longrightarrow\mathbb{N}$ is the function defined by
$f(0)=1$ and $f(r)=\dd\prod_{s=1}^r (p^s-1)$ for all
$r\in\mathbb{N}^*$.}
\bigskip

Obviously, explicit formulas for $C(\mathbb{Z}_p^n)$ and
$D(\mathbb{Z}_p^n)$ also follow from Theorems 1 and 2. By using a
computer algebra program, we are now able to calculate the first
terms of the chain $f_n=F(\mathbb{Z}_p^n)$, $n\in\mathbb{N}$,
namely:
\begin{description}
\item[\hspace{10mm}-] $f_0=1$;
\item[\hspace{10mm}-] $f_1=2$;
\item[\hspace{10mm}-] $f_2=2p+4$;
\item[\hspace{10mm}-] $f_3=2p^3+8p^2+8p+8$;
\item[\hspace{10mm}-] $f_4=2p^6+12p^5+24p^4+36p^3+36p^2+24p+16$.
\end{description}
\bigskip

Finally, we remark that the above $f_3$ is in fact the number
$a_{3,p}$ obtained by a direct computation in Corollary 10 of
\cite{5}.

\section{Proof of Theorem 3}

We observe first that every rooted chain of subgroups of
$\mathbb{Z}_p^n$ are of one of the following types:
$$G_1 \subset G_2 \subset ...\subset G_m=\mathbb{Z}_p^n\, \mbox{ with }\, G_1\neq 1\0(1)$$and
$$1 \subset G_2 \subset ...\subset
G_m=\mathbb{Z}_p^n\,.\0(2)$$It is clear that the numbers of chains
of types (1) and (2) are equal. So $$f_n=2x_n\,,\0(3)$$where $x_n$
denotes the number of chains of type (2). On the other hand, such
a chain is obtained by adding $\mathbb{Z}_p^n$ to the chain
$$1 \subset G_2 \subset ...\subset G_{m-1},$$where $G_{m-1}$ runs
over all subgroups of $\mathbb{Z}_p^n$. Moreover, $G_{m-1}$ is
also an ele\-mentary abelian $p$-group, say
$G_{m-1}\cong\mathbb{Z}_p^k$ with $0\leq k\leq n$. These show that
the chain $x_n$, $n\in\mathbb{N}$, satisfies the following
recurrence relation
$$x_n=\sum_{k=0}^{n-1} a_{n,p}(k)x_k\,,\0(4)$$which is more facile
than the recurrence relation founded by applying the
Inclusion-Exclusion Principle in Theorem 9 of \cite{5}.
\bigskip

Next we prove that the solution of (4) is given by
$$x_n=1{+}\sum_{k=1}^{n-1}\sum_{1\leq i_1<i_2<...<i_k\leq n-1} a_{n,p}(i_k)a_{i_k,p}(i_{k-1})\cdots
a_{i_2,p}(i_1)\,.\0(5)$$We will proceed by induction on $n$.
Clearly, (5) is trivial for $n=1$. Assume that it holds for all
$k<n$. One obtains
$$\hspace{-20mm}x_n=\sum_{k=0}^{n-1} a_{n,p}(k)x_k=1+\sum_{k=1}^{n-1} a_{n,p}(k)x_k=$$
$$\hspace{-3mm}=1{+}\sum_{k=1}^{n-1} a_{n,p}(k)\left(1+\sum_{r=1}^{k-1}\sum_{1\leq i_1<i_2<...<i_r\leq k-1} a_{k,p}(i_r)a_{i_r,p}(i_{r-1})\cdots a_{i_2,p}(i_1)\right)=$$
$$=1{+}\sum_{k=1}^{n-1} a_{n,p}(k)+\sum_{k=1}^{n-1}a_{n,p}(k)\sum_{r=1}^{k-1}\sum_{1\leq i_1<i_2<...<i_r\leq k-1} \hspace{-1mm}a_{k,p}(i_r)a_{i_r,p}(i_{r-1})\cdots a_{i_2,p}(i_1)=$$
$$=1{+}\sum_{k=1}^{n-1} a_{n,p}(k)+\sum_{k=1}^{n-1}a_{n,p}(k)\sum_{r=1}^{n-2}\sum_{1\leq i_1<i_2<...<i_r\leq k-1} \hspace{-1mm}a_{k,p}(i_r)a_{i_r,p}(i_{r-1})\cdots a_{i_2,p}(i_1)=$$
$$=1{+}\sum_{k=1}^{n-1} a_{n,p}(k)+\sum_{k=1}^{n-1}a_{n,p}(k)\sum_{r=2}^{n-1}\sum_{1\leq i_1<i_2<...<i_{r-1}\leq k-1} \hspace{-1mm}a_{k,p}(i_{r-1})a_{i_{r-1},p}(i_{r-2})\cdots a_{i_2,p}(i_1)=$$
$$=1{+}\hspace{-3mm}\sum_{1\leq i_1\leq n-1} a_{n,p}(i_1)+\sum_{r=2}^{n-1}\sum_{k=1}^{n-1}a_{n,p}(k)\hspace{-5mm}\sum_{1\leq i_1<i_2<...<i_{r-1}\leq k-1} \hspace{-2mm}a_{k,p}(i_{r-1})a_{i_{r-1},p}(i_{r-2})\cdots a_{i_2,p}(i_1)=$$
$$\hspace{-7mm}=1{+}\hspace{-3mm}\sum_{1\leq i_1\leq n-1} a_{n,p}(i_1)+\sum_{r=2}^{n-1}\sum_{1\leq i_1<i_2<...<i_r\leq n-1} a_{n,p}(i_r)a_{i_r,p}(i_{r-1})\cdots a_{i_2,p}(i_1)=$$
$$\hspace{14mm}=1+\sum_{r=1}^{n-1}\sum_{1\leq i_1<i_2<...<i_r\leq n-1} a_{n,p}(i_r)a_{i_r,p}(i_{r-1})\cdots a_{i_2,p}(i_1)\,,$$as desired.

Since by Theorem 2
$$a_{n,p}(k)=\frac{(p^n-1)\cdots (p-1)}{(p^k-1)\cdots (p-1)(p^{n-k}-1)\cdots (p-1)}=\frac{f(n)}{f(k)f(n-k)}\,,\forall\hspace{1mm} 0\leq k\leq
n\,,$$the equalities (3) and (5) imply that
$$f_n=2+2f(n)\sum_{k=1}^{n-1}\sum_{1\leq i_1<i_2<...<i_k\leq n-1}\frac{1}{f(n{-}i_k)f(i_k{-}i_{k-1})\cdots
f(i_2{-}i_1)f(i_1)}\,,$$completing the proof. $\scriptstyle\Box$

\vspace*{5ex}\small

\hfill
\begin{minipage}[t]{4cm}
Marius T\u arn\u auceanu \\
Faculty of  Mathematics \\
``Al.I. Cuza'' University \\
Ia\c si, Romania \\
e-mail: {\tt tarnauc@uaic.ro}
\end{minipage}

\end{document}